\documentclass[11pt]{article}

\usepackage{amscd,amsmath, amssymb, fancyhdr, graphics}

\numberwithin{equation}{section}

\newcommand{\version}{version 2.0,\ \  Oct 15, 2010}

\def\eqref#1{(\ref{#1})}
\newcommand{\goth}{\mathfrak}
\newcommand{\g}{{\mathfrak g}}
\newcommand{\arrow}{{\:\longrightarrow\:}}

\def\C{{\Bbb C}}
\newcommand{\R}{{\Bbb R}}

\newcommand{\6}{\partial}
\def\1{\sqrt{-1}\:}
\newcommand{\restrict}[1]{{\left|_{{\phantom{|}\!\!}_{#1}}\right.}}
\newcommand{\cntrct}{\rfloor}                

\renewcommand{\bar}{\overline}
\renewcommand{\phi}{\varphi}
\renewcommand{\epsilon}{\varepsilon}

\renewcommand{\leq}{\leqslant}

\newcommand{\hor}{{\operatorname{\rm hor}}}
\newcommand{\ver}{{\operatorname{\rm vert}}}

\newcommand{\End}{\operatorname{End}}
\newcommand{\Tot}{\operatorname{Tot}}
\newcommand{\Id}{\operatorname{Id}}

\newcommand{\Hom}{\operatorname{Hom}}
\newcommand{\Sym}{\operatorname{Sym}}

\newcommand{\Lie}{\operatorname{Lie}}

\newcommand{\Tw}{\operatorname{Tw}}

\newcounter{Mycounter}[section]
\newcounter{lemma}[section]
\renewcommand{\thelemma}{{Lemma \thesection.\arabic{lemma}}}
\newcommand{\lemma}{%
    \setcounter{lemma}{\value{Mycounter}}
    \refstepcounter{lemma}
    \stepcounter{Mycounter}
    {\noindent \bf \thelemma:\ }}

\newcounter{claim}[section]
\renewcommand{\theclaim}{{Claim \thesection.\arabic{claim}}}
\newcommand{\claim}{%
    \setcounter{claim}{\value{Mycounter}}
    \refstepcounter{claim}
    \stepcounter{Mycounter}
    {\noindent \bf \theclaim:\ }}

\newcounter{sublemma}[section]
\newcounter{corollary}[section]
\newcounter{theorem}[section]
\renewcommand{\thetheorem}{{Theorem \thesection.\arabic{theorem}}}
\newcommand{\theorem}{%
    \setcounter{theorem}{\value{Mycounter}}
    \refstepcounter{theorem}
    \stepcounter{Mycounter}
    {\noindent \bf \thetheorem:\ }}

\newcounter{conjecture}[section]
\newcounter{proposition}[section]
\renewcommand{\theproposition}
      {{Proposition \thesection.\arabic{proposition}}}
\newcommand{\proposition}{%
    \setcounter{proposition}{\value{Mycounter}}
    \refstepcounter{proposition}
    \stepcounter{Mycounter}
    {\noindent \bf \theproposition:\ }}

\newcounter{definition}[section]
\renewcommand{\thedefinition}
      {{Definition~\thesection.\arabic{definition}}}
\newcommand{\definition}{%
    \setcounter{definition}{\value{Mycounter}}
    \refstepcounter{definition}
    \stepcounter{Mycounter}
    {\noindent \bf \thedefinition:\ }}

\newcounter{example}[section]
\renewcommand{\theexample}{{Example \thesection.\arabic{example}}}
\newcommand{\example}{%
    \setcounter{example}{\value{Mycounter}}
    \refstepcounter{example}
    \stepcounter{Mycounter}
    {\noindent \bf \theexample:\ }}

\newcounter{remark}[section]
\renewcommand{\theremark}{{Remark \thesection.\arabic{remark}}}
\newcommand{\remark}{%
    \setcounter{remark}{\value{Mycounter}}
    \refstepcounter{remark}
    \stepcounter{Mycounter}
    {\noindent \bf \theremark:\ }}

\newcounter{problem}[section]
\newcounter{question}[section]
\makeatletter

\@addtoreset{equation}{section} \@addtoreset{footnote}{section}
\makeatother

\def\blacksquare{\hbox{\vrule width 5pt height 5pt depth 0pt}}
\def\endproof{\blacksquare}

\begin{document}
\begin{center}
{\LARGE\bf
A CR twistor space of a $G_2$-manifold\\[3mm]
}

 Misha
Verbitsky\footnote{The work is partially supported by 
RFBR grant 09-01-00242-a, Science Foundation of the 
SU-HSE award No. 09-09-0009 and by RFBR grant
10-01-93113-NCNIL-a.


}

\end{center}

{\small \hspace{0.15\linewidth}
\begin{minipage}[t]{0.7\linewidth}
{\bf Abstract} \\ Let $M$ be a $G_2$-manifold.
We consider an almost CR-structure on  the sphere
bundle of unit tangent vectors on $M$, called
{\bf the CR twistor space}. This CR-structure
is integrable if and only if $M$ is a holonomy
$G_2$ manifold. We interpret $G_2$-instanton
bundles as CR-holomorphic bundles on its
twistor space.
\end{minipage}
}

\tableofcontents


\section{Introduction}


\subsection{A CR twistor space of a $G_2$-manifold}

This note is inspired by the Claude LeBrun's paper
\cite{_Le_Brun:twistor_}. In his work, LeBrun 
constructed a twistor space for a 3-manifold, which happens
to be a CR-manifold of real dimension 5. The geometry of $G_2$-manifolds
is remarkably similar to the geometry of 3-manifolds, and
it is not surprising that an analogue of LeBrun's construction
can be obtained.

\hfill

\definition
Let $M$ be a smooth manifold, $B\subset TM$ a subbundle
in its tangent bundle, and $I\in \End B$ its automorphism,
$I^2=-\Id_B$. Consider the (1,0) and (0,1)-bundles
$B^{1,0}, B^{0,1}\subset B \otimes \C$, which are the
eigenspaces of $I$ corresponding to the eigenvalues
$\1$ and $-\1$. The sub-bundle $B^{1,0}\subset TM \otimes \C$ 
is called {\bf a CR-structure on $M$} if it is {\em involutive},
that is, satisfies $[B^{1,0}, B^{1,0}]\subset B^{1,0}$.

\hfill

\example
If $B=TM$, CR-structures are the same as complex
structures. For any codimension 1 real submanifold $M$
in a complex manifold $(X,I)$, the intersection $B:=TX \cap I(TX)$
is a complex subbundle of codimension 1 in $TX$, and
the restriction $I \restrict B$ defines a CR-structure.

\hfill

$G_2$-manifolds originally appeared
in Berger's classification of holonomy
(\cite{_Berger:list_}, \cite{_Besse:Einst_Manifo_}).
The first examples of $G_2$-manifolds were
obtained by R. Bryant and S. Salamon (\cite{_Bryant_Salamon_}).
The compact examples of $G_2$-manifolds
were constructed by D. Joyce (\cite{_Joyce_G2_}, 
\cite{_Joyce_Book_}). 
In this introduction we follow the approach
to $G_2$-geometry which is due to
N. Hitchin (see \cite{_Hitchin:3-forms_}).

\hfill

\definition
Let $\rho\in \Lambda^2 \R^7$ be a 3-form on $\R^7$.
We say that $\rho$ is {\bf non-degenerate} if the 
dimension of its stabilizer is maximal:
\[
\dim St_{GL(7)}\rho = \dim GL(7) - \dim \Lambda^3(\R^7) =
49-35 =14.
\]
In this case, $St(\rho)$ is one of two real forms 
of a 14-dimensional Lie group $G_2(\C)$. 
We say that $\rho$ is {\bf non-split} if it satisfies
$St(\rho|_x)\cong G_2$, where $G_2$ denotes the
compact real form of $G_2(\C)$.
{\bf A $G_2$-structure} on a 7-manifold
is a 3-form $\rho \in \Lambda^3(M)$, which is
non-degenerate and non-split at each point $x\in M$. We shall always
consider a $G_2$-manifold as a Riemannian
manifold, with the Riemannian structure induced
by the $G_2$-structure as indicated below.

\hfill

\remark \label{_metric_volume_G_2_Remark_}
Such a form defines
a $\Lambda^7 M$-valued metric on $M$:
\begin{equation}\label{_metric_vol_valued_Equation_}
g(x,y) = (\rho\cntrct x)\wedge (\rho \cntrct y) \wedge \rho
\end{equation}
(we denote by $\rho \cntrct x$ the contraction of $\rho$
with a vector field $x$). The Riemannian volume form
associated with this metric gives a section of 
$\Lambda^7 M\otimes (\Lambda^7 M)^{7/2}$. Squaring and
taking the 9-th degree root, we obtain a trivialization of the 
volume. Then \eqref{_metric_vol_valued_Equation_} defines
a metric $g$ on $M$, by construction $G_2$-invariant.

\hfill

\definition 
An  $G_2$-structure is called {\bf an integrable $G_2$-\-struc\-ture},
if $\rho$ is preserved by the corresponding Levi-Civita connection.
An integrable 
$G_2$-manifold is a manifold equipped with an integrable $G_2$-structure.
Holonomy group of such a manifold clearly lies in $G_2$;
for this reason, the integrable $G_2$-manifolds are often called
{\bf holonomy $G_2$-manifolds}. 

\hfill

\remark
In the literature, ``the $G_2$-manifold'' often means
a ``holonomy $G_2$-manifold'', and ``$G_2$-structure''
``an integrable $G_2$-structure''. A $G_2$-structure
which is not necessarily integrable is called 
``an almost $G_2$-structure'', taking analogy
from almost complex structures.

\hfill

\remark
As shown in \cite{_Fernandez_Gray:G_2_},
integrability of a $G_2$-structure induced
by a 3-form $\rho$ is equivalent to $d\rho = d(*\rho)=0$.
For this reason the 4-form $*\rho$ is called
{\bf a fundamental 4-form of a $G_2$-manifold},
and $\rho$ {\bf the fundamental 3-form}.

\hfill

\remark\label{_SU(3)_Remark_}
Let $V=\R^7$ be a 7-dimensional real space equipped
with a non-degenerate 3-form $\rho$ with $St_{GL(7)}(\rho)=G_2$.
As in \ref{_metric_volume_G_2_Remark_}, one can easily
see that $V$ has a natural $G_2$-invariant metric. 
For each vector $x\in V$, $|x|=1$, its stabilizer 
$St_{G_2}(x)$ in $G_2$ is isomorphic to $SU(3)$.
Indeed, the orthogonal complement $x^\bot$ is
equipped with a symplectic form $\rho \cntrct x$,
which gives a complex structure $g^{-1} \circ (\rho \cntrct x)$
as usual. This gives an embedding $St_{G_2}(x)\hookrightarrow U(3)$.
Since the space of such $x$ is $S^6$, and
the action of $G_2$ in $S^6$ is transitive, one has
$\dim St_{G_2}(x)= \dim G_2 - \dim S^6 =8 = \dim U(3)-1$.
To see that $St_{G_2}(x)= SU(3)\subset U(3)$ and not some
other codimension 1 subgroup, one should notice that
$St_{G_2}(x)$ preserves two 3-forms $\rho \restrict{x^\bot}$
and $\rho^*\cntrct x \restrict{x^\bot}$, where $\rho^* = *\rho$
is the fundamental 4-form of $V$. A simple linear-algebraic
calculation implies that $\rho \restrict{x^\bot}+\1 
\rho^*\cntrct x \restrict{x^\bot}$ is a holomorphic
volume form on $x^\bot$, which is clearly preserved by 
$St_{G_2}(x)$. Therefore, the natural embedding
$St_{G_2}(x)\hookrightarrow U(3)$ lands $St_{G_2}(x)$
to $SU(3)$. Using the dimension count $\dim St_{G_2}(x)=\dim SU(3)$
(see above), we show that the embedding
$St_{G_2}(x)\hookrightarrow SU(3)$ is also surjective.

\hfill

Let now $M$ be an almost $G_2$-manifold. 
From \ref{_SU(3)_Remark_} it follows that with every
vector $x\in TM$, $|x|=1$, one can associate a
complex Hermitian structure on its orthogonal
complement $x^\bot$. The easiest way to define 
this structure is to notice that $x^\bot$
is equipped with a symplectic structure 
$\rho\cntrct x$ and a metric $g\restrict{x^\bot}$,
which can be considered as a real and imaginary
parts of a complex-valued semilinear Hermitian 
product. Then the complex structure is obtained as
usual, as $I:=(\rho\cntrct x)\circ g^{-1}$.

\hfill

\definition\label{_B^1,0_Definition_}
Consider now the unit sphere bundle $S^6M$
over $M$, with the fiber $S^6$, and let $T_\hor S^6M \subset TS^6 M$
be the horizontal sub-bundle corresponding to the
Levi-Civita connection. This sub-bundle has a natural
section $\theta$; at each point $(x, m)\in S^6 M$,
$m\in M$, $x\in T_m M$, $|x|=1$, we take
$\theta\restrict{(x, m)}= x$, using the standard 
identification $T_\hor S^6 M\restrict{(x, m)}= T_m M$.
Denote by $B\subset T_\hor S^6 M$ the orthogonal complement to $\theta$
in $T_\hor S^6 M$. Since at each point $(x, m)\in S^6 M$,
the restriction $B \restrict{(x, m)}$ is identified
with $x^\bot \subset T_m M$, this bundle is equipped
with a natural complex structure, that is, an operator
$I \in \End B$, $I^2= -\Id_B$. 

\hfill

The main result of the present paper is the
following theorem.

\hfill

\theorem\label{_CR_main_Theorem_}
Let $M$ be an almost $G_2$-manifold, $S^6M\subset TM$
its unit sphere bundle, and $B\subset T S^6 M$
a sub-bundle of its tangent bundle constucted
above, and equipped with the complex structure $I$ as above
Then $B^{0,1}\subset B\otimes \C \subset T S^6 M\otimes \C$
is involutive if ang only if $M$ is a holonomy $G_2$-manifold.

\hfill

We prove \ref{_CR_main_Theorem_} in 
Section \ref{_integra_proof_Section_}.

\hfill

\definition\label{_twistor_Definition_}
Let $M$ be a holonomy $G_2$-manifold,
and \[ \Tw(M):=(S^6 M, B, I)\] the CR-manifold constructed 
in \ref{_CR_main_Theorem_}. Then $\Tw(M)$ is called
{\bf a CR-twistor space of $M$}.

\subsection{Applications of twistor geometry}

$G_2$ instanton bundles were introduced in 
\cite{_Donaldson_Thomas_}, and much studied since then.
This notion is a special case of a more general
notion of an instanton on a calibrated manifold,
which is already well developed. Many estimates 
known for 4-dimensional manifolds
(such as Uhlenbeck's compactness theorem)
can be generalized to the calibrated case 
(\cite{_Tian:Calibrated_}, \cite{_Tian_Tao_}).

Recently, $G_2$-instantons became a focus of
much activity because of attempts to construct
a higher-dimensional topological quantum field theory,
associated with $G_2$ and 3-dimensional Calabi-Yau
manifolds (\cite{_Donaldson_Segal_}).

\hfill

\definition
Let $M$ be a $G_2$-manifold, and 
$\Lambda^2 M = \Lambda^2_7(M) \oplus \Lambda^2_{14}(M)$
the irreducible decomposition of the bundle of 2-forms 
$\Lambda^2(M)$ associated with the $G_2$-action
(Subsection \ref{_G_2_decompo_Subsection_}). 
A vector bundle $(B, \nabla)$ with connection
is called {\bf a $G_2$-instanton} if its curvature
lies in $\Lambda^2_{14}(M)\otimes \End(B)$.

\hfill

\remark
The instanton connection is the same as a self-dual
connection in the context of $G_2$-geometry. In particular,
it minimizes the curvature functional, in the same way
as the usual instantons minimize the curvature
functional. Therefore, the instanton connections
are solutions of the Yang-Mills equation
(see \cite{_Tian:Calibrated_}).

\hfill

\remark
A tangent bundle and all its tensor powers are obviously
$G_2$-instantons (see Step 1 in the proof of \ref{_CR_main_Theorem_}
in Section \ref{_integra_proof_Section_}). 

\hfill

The definition of CR-holomorphic bundles is a 
straighforward generalization of a usual differential-geometric
notion of a holomorphic bundle as a bundle equipped with
a Dolbeault  differential $\bar\6$ which satisisfies
$\bar\6^2=0$. 

\hfill

\definition
Let $(M, B, I\in \End B)$ be a CR-manifold,
and $B^{0,1}\subset B\otimes \C$ the $\1$-eigenspace of $B$.
Using the Cartan's formula, 
we define the CR Dolbeault differential 
\[ 
  \Lambda^k(B^{0,1})^* \stackrel {\bar\6_B}\arrow \Lambda^{k+1}(B^{0,1})^*
\]
as usual,
\begin{align*}
(\bar\6_B\alpha)(b_1, ... b_{k+1}):= &
\sum_i (-1)^i \Lie_{b_i}\alpha(b_1, ..., \check b_i,..., b_{k+1})\\
+ & \sum_{i<j} (-1)^{i+j+1} \alpha([b_i, b_j], b_1, ..., \check 
b_i, ..., \check b_j, ..., b_{k+1}).
\end{align*}
Let $E$ be a vector bundle, and $\bar\6_E:\; E \arrow E\otimes (B^{0,1})^*$
an operator which satisfies the Leibniz rule 
\[ 
  \bar\6_E(f\xi) = \xi \otimes \bar\6_Bf + f \bar\6_E\xi.
\]
We extend $\bar\6_E$ to 
\[
 E\otimes \Lambda^k(B^{0,1})^* \stackrel {\bar\6_E} \arrow 
E\otimes \Lambda^{k+1}(B^{0,1})^*
\]
using the same Leibniz formula,
\[
\bar\6_E(\xi \otimes \alpha) = 
\xi \otimes \bar\6_B\alpha + \bar\6_E\xi\otimes \alpha.
\]
for any $\xi \in E, \alpha \in \Lambda^k(B^{0,1})^*$.
A bundle $(E, \bar\6_E)$ is called {\bf CR-holomorphic}, if
$\bar\6_E^2=0$. In this case, $\bar\6_E$ is called
{\bf an operator of CR-holomorphic structure}.

\hfill

\definition
Let $(M, B, I\in \End B)$ be a CR-manifold, and
$(E, \nabla)$ a complex vector bundle on $M$ with connection.
Denote by \[ \Pi:\; \Lambda^i M \arrow  \Lambda^i(B^{0,1})^*\]
the restriction map, and let 
\[ \bar\6_E:=\nabla \circ\Pi:\; E \arrow E\otimes (B^{0,1})^*\]
be the $B^{0,1}$-part of the connection. We say that $(E,\nabla)$
is {\bf CR-holomorphic} if $\bar\6_E$ is an operator of
CR-holomorphic structure on $E$.

\hfill

\theorem
Let $M$ be a holonomy $G_2$-manifold, $(E, \nabla)$
a bundle with connection, and $\Tw(M)\stackrel\pi \arrow M$
its twistor space defined as in \ref{_twistor_Definition_}.
Then the following assertions are equivalent.
\begin{description}
\item[(i)] The pullback $(\pi^*E, \pi^*\nabla)$ is
a CR-holomorphic bundle on $\Tw(M)$
\item[(ii)] $(E, \nabla)$ is a $G_2$-instanton.
\end{description}

{\bf Proof:} Let $\nabla^2 R \in \Lambda^2M \otimes \End E$
be the curvature form of $E$, $\pi^*R$ the curvature
of $(\pi^*E, \pi^*\nabla)$, and
$\Pi(\pi^* R) \in \Lambda^2(B^{0,1})^* \otimes \End E$
its restriction to $B^{0,1}$. Clearly, $\Pi(\pi^* R)=0$
if and only if $(\pi^*E, \pi^*\nabla)$ is CR-holomorphic. 
However, $\Pi(\pi^* R)=0$ if and only if $\pi^* R$
is of Hodge type $(1,1)$ on $B$, and this is equivalent
to $R\in \Lambda^2_{14}\otimes \End E$, as follows from
\ref{_Lambda^2_14_via_Hodge_Proposition_}.
\endproof

\hfill

\remark
Following the same approach as in the paper \cite{_Lempert:CR_},
it is possible to deduce from \ref{_CR_main_Theorem_}
that the space of knots (non-parametrized loops) in
a holonomy $G_2$ manifold is formally Kaehler (see
the forthcoming paper \cite{_Verbitsky:Knot_}). For 3-manifolds
this theorem is due to J.-L. Brylinski (\cite{_Brylinski_}).


\section{Differential forms on $\Tot(\Lambda^k M)$}
\label{_forms_on_Tot_Lambda^k_Section_}

In this section, we extend the standard results about the
Hamiltonian 2-form on $T^*M$ to the total space $\Tot(\Lambda^k(M))$
of the bundle of $k$-forms. These generalizations are elementary,
but used further on in this paper.

\hfill

Let $M$ be a smooth manifold, and $X= \Tot(\Lambda^k M)$
a total space of a bundle of $k$-forms. On $X$, a pair
of forms is defined: a $k$-form $\Theta$, which is an 
analogue of the 1-form $\sum q_i dp_i$ on $T^*M$, 
and a $(k+1)$-form $\Xi=d\Theta$, which is an analogue
of the symplectic form.

Let $\pi:\; X \arrow M$ be a standard projection.
At a point $(\lambda, m)\in X$, with $m\in M$ and 
$\lambda\in \Lambda^k(M)\restrict m$, we define 
\begin{equation}\label{_Theta_in_coo_Equation_}
   \Theta(x_1, ..., x_k):= 
   \lambda(D\pi(x_1), D\pi (x_2), ..., D\pi(x_k)),
\end{equation}
where $D\pi:\; TX \arrow TM$ is the differential of $\pi$.

Let $\Xi:= d\Theta$. In local coordinates $p_1, ... p_n$ on
$M$, 
\[ 
q_{i_1, ..., i_k}:= 
   dp_{i_1}\wedge dp_{i_2}\wedge ... \wedge dp_{i_k}
\]
on the fibers of $\pi$, $\Xi$ can be written as as sum
\begin{equation}\label{_Xi_in_coo_Equation_}
\Xi = 
\sum_{i_1 < i_2< ... i_k} dq_{i_1, ..., i_k}\wedge  dp_{i_1}\wedge dp_{i_2}\wedge ... \wedge dp_{i_k}
\end{equation}
(this is clear from the same argument as used to obtain
a more familiar relation $d(\sum q_i dp_i) = \sum dq_i \wedge dp_i$).
Further on, we shall need the following elementary lemma.

\hfill

\lemma\label{_Xi_on_hori_Lemma_}
Let $M$ be a manifold equipped with a torsion-free connection $\nabla$,
$X= \Tot(\Lambda^k M)$ the total space
of the bundle of $k$-forms, and $\Xi\in \Lambda^{k+1}(X)$
the fundamental $(k+1)$-form constructed above. 
Using the connection $\nabla$, we split the
tangent bundle onto horizontal and vertical components,
$TX = T_\hor X \oplus T_\ver X$. Then $\Xi\restrict{T_\hor X}=0$.

\hfill

{\bf Proof:} Choose coordinates $p_1, ..., p_n$
on $M$, and let $p_1, ... p_n$, $q_{i_1, ..., i_k}$,
$1 \leq i_1 < i_2 < ... < i_k \leq n$ be the coordinates
on $X$ constructed as above.  This coordinate system
induces a flat torsion-free connection $\nabla_0$ on $TM$ and
$\Lambda^k M$.
For an arbitrary
connection $\nabla$ on $TM$, one obtains a splitting
$TX = T_\ver \oplus T_\hor(\nabla)$.
Write $\nabla=\nabla_0 + A$, where $A\in \Lambda^1M \otimes \End(TM)$
is an $\End(TM)$-valued 1-form on $M$.  The difference between 
$T_\hor(\nabla)$ and $T_\hor(\nabla_0)$ can be expressed as a 
section of $\Hom(\pi^* TM, T_\ver X)$, at $(\lambda, m)$
giving $A_k(\lambda) \in \Lambda^1(M)\otimes T_\ver X$,
where $A_k\in \Lambda^1(M)\otimes \End(\Lambda^k M)$
is the connection form on $\Lambda^k(M)$ induced
from $\nabla$, $A_k = \nabla-\nabla_0$.
Denote $A_k(\lambda)$ by 
$E:\; \pi^*\Lambda^1M\arrow T_\ver X=\pi^* \Lambda^k M$. 
Using the formula for $\Xi$ written
in coordinates as in \eqref{_Xi_in_coo_Equation_},
we obtain 
\begin{equation}\label{_Xi_on_hor_Equation_}
\Xi(v_{i_1}', v_{i_2}', ...,
v_{i_{k+1}}')= \Xi(v_{i_1}+E(v_{i_1}), v_{i_2}+E(v_{i_2}), ...,
v_{i_{k+1}}+E(v_{i_{k+1}})),
\end{equation}
where $v_i:= \frac d {dp_i}$ are coordinate vector fields
on $X$, and $v_i'= v_i + E(v_i)$ the corresponding
sections of $T_\hor(\nabla)$.
As follows from \eqref{_Xi_in_coo_Equation_},
the latter expression gives
\begin{equation}\label{_Xi_on_hor_permuta_Equation_}
\Xi(v_{i_1}', v_{i_2}', ...,
v_{i_{k+1}}') = \frac 1 {(k+1)!}\sum
(-1)^{|\sigma|}E(v_{\sigma(i_1)})\bigg(v_{\sigma(i_2)},v_{\sigma(i_3)},
...,  v_{\sigma(i_{k+1})}\bigg)
\end{equation}
where the sum is taken over all permutations $\sigma \in S_{k+1}$.
The operator $E:\; \pi^*\Lambda^1M\arrow \pi^* \Lambda^k M$
can be expressed in terms of the connection form
 $A\in \Lambda^1M \otimes \End(TM)$ as follows:
at $(\lambda, m) \in X$ one has
\begin{equation}\label{_E_formula_conne_Equation_}
\begin{aligned}
E\restrict{(\lambda, m)}(v)(v_{i_1}, ..., v_{i_k})\\ =&
\lambda(A_v(v_{i_1}), v_{i_2}, ..., v_{i_k}) +
\lambda(v_{i_1}, A_v(v_{i_2}),  ..., v_{i_k}) \\ + &
... + \lambda(v_{i_1},v_{i_2}, ...,  A_v(v_{i_k})),
\end{aligned}
\end{equation}
where $A_v$ denotes $A(v)$.
Comparing \eqref{_E_formula_conne_Equation_} and
\eqref{_Xi_on_hor_permuta_Equation_}, we obtain
\begin{multline}\label{_Xi_on_hor_via_connec_form_Equation_}
\Xi(v_{i_1}', v_{i_2}', ...,
v_{i_{k+1}}')\restrict{(\lambda, m)} \\ = \frac 1 {(k+1)!}\sum
(-1)^{|\sigma|}\lambda\bigg(A_{v_{\sigma(i_1)}}(v_{\sigma(i_2)}),v_{\sigma(i_3)},
...,  v_{\sigma(i_{k+1})}\bigg).
\end{multline}
Since $\nabla$ and $\nabla_0$ are torsion-free, one
has $A_{v_{\sigma(i_1)}}(v_{\sigma(i_2)})=A_{v_{\sigma(i_2)}}(v_{\sigma(i_1)})$.
Therefore, the alternating sum 
\eqref{_Xi_on_hor_via_connec_form_Equation_}
vanishes. \endproof

\hfill

\remark 
This argument repeats a more familiar argument showing that
the restriction of the Hamiltonian symplectic form 
$\omega=\sum_i dp_i\wedge dq_i$
on $T^*M$ to a horizontal space of a 
torsion-free connection is always zero.


\section[The $G_2$-action on $\Lambda^2(\R^7)$ and
$SU(3)$-action  on $\Lambda^2(\R^6)$]{The $G_2$-action on $\Lambda^2(\R^7)$ and
$SU(3)$-action \\on $\Lambda^2(\R^6)$}


\subsection{Octonion algebra and quaternions}

Let $V=\R^7$ be a 7-dimensional space equipped
with a non-degenerate, non-split 3-form $\rho$ 
inducing a $G_2$-action on $V$. Then $V$ is equipped with
the vector product, defined as follows: 
$x\star y = \rho(x, y, \cdot)^\sharp$. Here
$\rho(x, y, \cdot)$ is a 1-form obtained by
contraction, and $\rho(x, y, \cdot)^\sharp$
its dual vector field. It is not hard to see
that $(V, \star)$ becomes isomorphic to the
imaginary part of the octonion algebra,
with $\star$ corresponding to half of the commutant.
In fact, this is one of a definitions of
an octonion algebra. The whole octonion algebra
is defined ${\Bbb O}:=V \oplus \R$,
with the product given by 
\[ (x, t) (y, t')= 
   (ty + t'x + x\star y, g(x,y)+tt')
\]
Here, $x, y$ and $ty + t'x + x\star y$
are vectors in $V$, and $t, t', g(x,y)+tt'\in \R$.

Given two non-collinear vectors in $V$, they generate
a quaternion subalgebra in octonions. When these two vectors
satisfy $|v|=|v'|=1$, $v\bot v'$, the standard basis $I, J, K$
in imaginary quaternions can be given by a triple 
$v, v', v\star v'\in V$.

A 3-dimensional subspace $A \subset V$ is called
{\bf associative} if it is closed under the vector
product. The set of associative subspaces is in
bijective correspondence with the set of quaternionic
subalgebras in octonions.

The following well-known lemma is used further on 
in this section.

\hfill

\lemma\label{_restri_hermi_hol_symple_Lemma_}
Let $(V, \rho)$ be a  7-dimensional space 
with a $G_2$-structure, $v\in V$ a unit
vector, and $W= v^\bot\subset V$
its orthogonal complement, equipped with
a Hermitian structure as explained in
\ref{_SU(3)_Remark_}. Consider another 
unit vector $v_1\in V$, $v_1 \bot v$,
and let $W_1:=v_1^\bot$, with its
own complex structure. Denote
by $\omega_W\in \Lambda^2 W$ the 
Hermitian form on $W$. Then the restriction
$\omega_W \restrict {W_1}$
is of Hodge type $(2,0)+(0, 2)$.

\hfill

{\bf Proof:} Let $A\subset V$ be a 3-dimensional
associative space generated by $v$ and $v_1$,
and $H= A^\bot$. The intersection $W\cap W_1$
is a 5-dimensional space generated by $H$ and
$v\star v_1$, where $\star$ denotes the
vector product. The restriction $\omega_W\restrict{W_1}$
vanishes on $v\star v_1$, because the complex structure
on $W$ maps $v\star v_1$ to $-v_1$. Therefore,
$\omega_W\restrict{W_1}$ is an image of
$\omega_W\restrict{H}$ under the standard orthogonal
embedding $\Lambda^2 H \arrow \Lambda^2W_1$.

Three complex structures $v$, $v_1$, $v\star v_1$
define a quaternionic structure on $H$, and
the corresponding three Hermitian forms
give a hyperkaehler triple of symplectic
structures on $H$. Since $\omega_W\restrict H$
is a real part of the standard $(2,0)$-form
on $H$ associated with $v_1$, it is
a $(2,0)+(0,2)$-form on $W_1$. \endproof

\subsection{$G_2$-action and $SU(3)$-action}
\label{_G_2_decompo_Subsection_}

Let $V=\R^7$ be a 7-dimensional space equipped
with a non-degenerate, non-split 3-form $\rho$ 
inducing a $G_2$-action on $V$ as above.
The space $\Lambda^2 V$ is a 21-dimensional
representation of $G_2$. The Lie algebra
$\g_2$ can be considered as a subspace in $\Lambda^2 V$,
because one has an embedding 
$\g_2 \hookrightarrow \goth{so}(V) =\Lambda^2 V$.
This gives a 14-dimensional $G_2$-invariant subspace
in $\Lambda^2 V$. 
There is also a 7-dimensional subspace given by
an embedding $V \arrow \Lambda^2 V$, $v \mapsto \rho \cntrct  v$.
Since $\Lambda^2 V$ is 21-dimensional, this gives a
decomposition $\Lambda^2 V = \Lambda^2_7 V\oplus \Lambda^2_{14} V$
of $\Lambda^2 V$ onto irreducible 7-dimensional and 
14-dimensional $G_2$-sub\-rep\-re\-sen\-ta\-tions.

Further on in this paper, we shall need the
following linear-algebric result.

\hfill

\proposition\label{_Lambda^2_14_via_Hodge_Proposition_}
Let $V=\R^7$ be a 7-dimensional space equipped
with a $G_2$-action, and $\Lambda^2 V = \Lambda^2_7 V\oplus \Lambda^2_{14} V$
the irreducible decomposition obtained as above. 
For any $v\in V$, consider the Hermitian structure
on its orthogonal complement $v^\bot$ (\ref{_SU(3)_Remark_}).
Let $\alpha \in \Lambda^2 V$ be a 2-form. Then
the following conditions are equivalent.
\begin{description}
\item[(i)] $\alpha$ lies in $\Lambda^2_{14} V$.
\item[(ii)] For any non-zero $v\in V$, the restriction
of $\alpha$ to $v^\bot$ is of type $(1,1)$
with respect to the complex structure on 
$v^\bot$, and orthogonal to the Hermitian form.
\item[(iii)] For any non-zero $v\in V$, the restriction
of $\alpha$ to $v^\bot$ is of type $(1,1)$
with respect to the complex structure on 
$v^\bot$.
\end{description}

{\bf Proof:} 
Denote by $W\subset V$ the orthogonal complement
to $v$, with the natural $SU(3)$-structure
constructed as in \ref{_SU(3)_Remark_}. 
After the standard identification between $\goth{so}(W)$
and $\Lambda^2 W$, the algebraic condition (ii)
is translated to $\alpha\restrict W \in \goth{su}(W)$.
The restriction map $\Lambda^2 V \arrow \Lambda^2 W$
is dual to the standard embedding 
$\goth{so}(W) \hookrightarrow \goth{so}(V)$.
Under this embedding, $\goth{su}(W)$ is mapped to
$\g_2$ (\ref{_SU(3)_Remark_}), hence 
the restriction $\Lambda^2_{14}(V)\restrict W$
lies in $\goth{su}(W)$.  We proved the
implication (i) $\Rightarrow$ (ii).

To obtain the converse implication, take some 
$\alpha\in \Lambda^2 V$ satisfying assumptions
of (ii); since $V$ is odd-dimensional,
$\alpha$ has an annulator $R\subset V$, which is at least 
one-dimensional. Let $v\in R$ be a non-zero vector.
Then $\alpha$, considered as an element
of $\goth{so}(V)$, preserves $W= v^\bot$
and acts trivially on $v$, hence it is an
image of $\alpha_0\in \Lambda^2(W)$
under the standard embedding 
$\goth{so}(W) \hookrightarrow \goth{so}(V)$.
By our assumptions, $\alpha_0 \in \goth{su}(W)$, 
hence $\alpha$ belongs to an image
of the standard embedding 
$\goth{su}(W) \hookrightarrow \g_2=\Lambda^2_{14}(V)$.

To finish the proof of
\ref{_Lambda^2_14_via_Hodge_Proposition_},
it suffices to prove that (iii) implies (i).
Let $\alpha \in \Lambda^2 V$ be a
2-form satisfying (iii). Consider
a codimension 1 subspace $W\subset V$
constructed as above, with $\alpha$
in the image of the standard orthogonal
embedding $\Lambda^2 W \hookrightarrow \Lambda^2 V$.
Consider a decomposition $\alpha = \alpha_0 + t\omega_W$
of $\alpha\in \goth{u}(W) = \Lambda^{1,1}(W)$ onto 
its traceless part $\alpha_0 \in \goth{su}(W)$
and the part $t\omega_W$ proportional to the Hermitian
form. Since $\alpha_0 \in \goth{su}(W)$, this
form satisfies $\alpha_0 \in \Lambda^2_{14}(V)$,
because $\goth{su}(W)$ lies in $\g_2$.
To prove that $\alpha \in \Lambda^2_{14}(V)$
it would suffice to show that $t\omega_W=0$.
However, since $\alpha_0$ lies in $\Lambda^2_{14}(V)$,
it satisfies the asusmptions of
\ref{_Lambda^2_14_via_Hodge_Proposition_}
(ii), hence satisfies \ref{_Lambda^2_14_via_Hodge_Proposition_}
(iii). Therefore, the same is true for $t\omega_W= \alpha-\alpha_0$:
the restriction of $t\omega_W$ to any $W'= (v')^\bot$ of codimension 1 in $V$
is of type $(1,1)$ with respect to the standard complex
structure on $W'$. This implies $t=0$, because,  as follows from 
\ref{_restri_hermi_hol_symple_Lemma_}, the restriction
of $\omega_W$ to another 6-dimensional subspace
$W'$ is of type $(2,0)+(0,2)$.
\endproof


\section{Integrability of the twistor CR-structure}
\label{_integra_proof_Section_}


The following claim is well-known from the standard
textbooks on differential geometry.

\hfill

\claim\label{_commu_conne_Claim_}
(\cite{_Kob_Nomizu_})
Let $M$ be a Riemannian $n$-manifold, $S^{n-1} M$
the sphere bundle consisting of all unit tangent vectors,
and $T_\hor S^{n-1}(M)\subset T S^{n-1}(M)$ the
horizontal tangent bundle associated with the
Levi-Civita connection. Denote by $R\in \Sym^2 \Lambda^2 M$
the curvature of the Levi-Civita connection on $M$.
Given a 2-form $\alpha \in \Lambda^2 M$, consider
$\alpha$ as an element in $\goth{so}(TM)=\Lambda^2 M$,
acting on $S^{n-1}(M)$ by infinitesimal automorphisms,
and let $\alpha^\star\in T_\ver S^{n-1} M$ be the corresponding
vector field. Let $X, Y \in T_\hor S^{n-1}(M)$
be horizontal vector fields. Then the vertical
component of $[X,Y]$ is equal to $R(X,Y)^\star$.
\endproof

\hfill

We are going to prove \ref{_CR_main_Theorem_},
which says that the bundle 
$B^{1,0}\subset T_\hor X\otimes \C$
on a $G_2$-manifold $M$ is involutive, where $X= S^6 M$.
The strategy of the proof is the following.
We consider the Frobenius bracket
$[B^{1,0}, B^{1,0}] \stackrel\Psi \arrow TX\otimes \C / B^{1,0}$.
First, we project $\Psi(X, Y)$, $X, Y \in B^{1,0}$
to $T_\ver X$ along the horizontal component
$T_\hor X \otimes \C \supset B^{1,0}$, showing that this
projection vanishes. This implies that $\Psi$ 
is a map from $\Lambda^2 B^{1,0}$ to 
$T_\hor X \otimes \C/B^{1,0}$:
\[ 
[B^{1,0}, B^{1,0}] \stackrel\Psi \arrow  T_\hor X\otimes \C / B^{1,0}
\]
Then we prove that the Frobenius bracket
$[\theta^\bot, \theta^\bot] \arrow TX/\theta^\bot$
vanishes at $B\subset TX$, where $\theta\in TX$ is a
natural section of $T_\hor X$ constructed as in
\ref{_B^1,0_Definition_}, and $\theta^\bot$
is its orthogonal complement. Since 
$B= T_\hor X \cap \theta^\bot$, this 
implies that the Frobenius bracket
$[B^{1,0}, B^{1,0}] \stackrel\Psi \arrow  T_\hor X\otimes \C / B^{1,0}$
lands in $\theta^\bot \cap T_\hor X = B$.
This brings us to the situation familiar
from complex (and almost complex) geometry,
where we have a bundle $B$ with an endomorphism
$I\in \End B$, $I^2=-\Id_B$, and the commutator map
$[B^{1,0}, B^{1,0}] \arrow B\otimes \C$
determines whether $I$ is integrable.
To prove that 
\begin{equation}\label{_commu_inte_Equation_}
[B^{1,0}, B^{1,0}] \subset B^{1,0},
\end{equation}
we construct a closed $(0,3)$-form $\bar\Omega$
which vanishes at $B^{1,0}$ and is non-\-de\-ge\-ne\-rate
at $B^{0,1}$, and use this form to prove 
\eqref{_commu_inte_Equation_}. 

\hfill

{\bf Proof of \ref{_CR_main_Theorem_}. Step 1.}
By \ref{_commu_conne_Claim_}, the commutator bracket 
\[
[T_\hor X, T_\hor X] \stackrel {\Psi_1}\arrow TX/T_\hor X
\]
is expressed through the curvature $R$ of the Levi-Civita
connection on $M$, as follows:
\[
\Psi_1(x,y)\restrict{(z,m)\in X}= R(x,y,z).
\]
Consider $R_z:= R(\cdot, \cdot, z)$
as a $TM$-valued 2-form on $M$. If $M$ is holonomy
$G_2$-manifold, this form lies in $\Lambda^2_{14}\otimes TM$,
because the curvature of a holonomy $G_2$ manifold lies in
$\Lambda^2_{14}\otimes \g_2$. Therefore, $R_z$ of type
$(1,1)$ on $B$, as follows from
\ref{_Lambda^2_14_via_Hodge_Proposition_}.
This implies that $[B^{1,0}, B^{1,0}]\subset T_\hor X\otimes \C$.
Conversely, if $[B^{1,0}, B^{1,0}]\subset T_\hor X\otimes \C$,
this means that $R_z$ is of type $(1,1)$
on every 6-dimensional subspace $W\subset T_m M$,
and \ref{_Lambda^2_14_via_Hodge_Proposition_}
implies that $R\in\Lambda^2_{14}\otimes \End(TM)$.

This already proves that integrability
of the CR structure $B^{1,0}\subset TX\otimes \C$
implies that $(M, \rho)$ is a holonomy $G_2$ manifold.
The converse implication takes more work.

\hfill

{\bf Step 2:} To study the bundle
$B = \theta^\bot \cap T_\hor X$, we identify
$X = S^6 TM$ with the total space of the
unit bundle $S^6 T^* M\subset \Tot(T^* M)$.
This is done using the Riemannian metric on
$TM$. On $\Tot(T^* M)$ there is a standard
1-form $\Theta$, dual to $\theta$
(see Section \ref{_forms_on_Tot_Lambda^k_Section_}). 
The zero-space of this form is identified with
$\theta^\bot$, and the
corresponding Frobenius bracket 
$[\theta^\bot,  \theta^\bot] \arrow TX /\theta^\bot$
is a 2-form which is proportional to the
restriction $d\Theta\restrict{\theta^\bot}$.
However, $d\Theta=\Xi$ is the standard
symplectic form on $T^*M$, and it
vanishes on the horizontal sub-bundle,
as follows from \ref{_Xi_on_hori_Lemma_}.
This implies that the bracket
$[\theta^\bot,  \theta^\bot] \arrow TX /\theta^\bot$
vanishes on $B \subset \theta^\bot$.
In Step 1, we have shown that 
$[B^{1,0}, B^{1,0}]\subset T_\hor X\otimes \C$;
the above argument implies that 
\[ 
   [B^{1,0}, B^{1,0}]\subset (T_\hor X\cap
   \theta^\bot)\otimes \C = B\otimes \C.
\]

{\bf Step 3:} The 6-dimensional bundle
$B \subset T_\hor M$ by construction comes 
equipped with a natural $SU(3)$-structure. 
The corresponding (3,0)-form 
can be written down as follows 
(\ref{_SU(3)_Remark_}):
\begin{equation}\label{_SU(3)_forms_Equation_}
  \bigg(\pi^* \rho + \1 (\pi^* \rho^*)
  \cntrct\theta\bigg)\restrict B
\end{equation}
Here, $\pi:\; X \arrow M$ is the standard
projection, $\rho$ and $\rho^*:=*\rho$ are the fundamental 3-form
and 4-form on $M$, and \eqref{_SU(3)_forms_Equation_}
is the standard equation expressing the standard $(3,0)$-form
on a codimension 1 subspace $W\subset V=\R^7$
in terms of the fundamental 3-form
and 4-form on $V$. Let 
$\Omega:= \pi^* \rho + \1 (\pi^* \rho^*)\cntrct\theta$.
In Step 3, we are going to 
prove that $(d\Omega)\restrict {T_\hor X}=0$. Since $\rho$
is closed, this would follow from 
$d\left((\pi^* \rho^*)  \cntrct\theta\right)\restrict {T_\hor X}=0$.
Let $Z:= \Tot(\Lambda^3 M)$, and $\phi:\; X \arrow Z$
be an embedding mapping $(v, m) \in X = S^6 M$
to $(\rho^* \cntrct v, m)$. Clearly, 
$\pi^*(\rho^*)=\phi^*(\Theta)$, where $\Theta$ is a 3-form
on $\Tot(\Lambda^3 M)$ defined in Section
\ref{_forms_on_Tot_Lambda^k_Section_}.
Then $d((\pi^* \rho^*)  \cntrct\theta)= \phi^*(\Xi)$.
However, the connection on $\Lambda^3 M$ is
compatible with the embedding 
$\Lambda^1 (M) \arrow\Lambda^3 M$,
hence the image $\phi(T_\hor X)$ lies
in the horizontal subspace $T_\hor Z\subset TZ$.
Therefore, to prove that 
\[
  d((\pi^* \rho^*)  \cntrct\theta)\restrict{T_\hor X} 
  = \phi^*(\Xi)\restrict{T_\hor X} =0
\]
it would suffice to show that $\Xi \restrict{T_\hor Z}=0$.
This follows from \ref{_Xi_on_hori_Lemma_}.

\hfill

{\bf Step 4:} 
Let $Z, T\in B^{0,1}$. Since $\Omega$ is of type
$(3,0)$ on $B$, one has $\Omega\cntrct Z= \Omega\cntrct T=0$.
Therefore, Cartan's formula gives
\begin{equation}\label{_Cartan_s_Equation_}
d\Omega(X,Y,Z,T)= \Omega(X,Y, [Z,T]).
\end{equation}
Taking a complex conjugate of  the equation
$[B^{1,0}, B^{1,0}]\subset B\otimes \C$ (Step 2),
we obtain that $[Z,T]\in B\otimes \C$.
Since $d\Omega\restrict B=0$, \eqref{_Cartan_s_Equation_}
implies that the 2-form $\Omega([Z,T], \cdot, \cdot)$
is identically zero on $B$. However, $\Omega$ is a non-degenerate
$(3,0)$-form on $B$, hence for any $\xi \in B\otimes \C$, the vanishing
$\Omega(\xi, \cdot, \cdot)=0$
implies that $\xi \in B^{0,1}$. We proved that
$B^{0,1}$ is involutive. The same is true
for $B^{1,0}$ by complex conjugation.
\endproof

\hfill

\hfill

\small{

\noindent {\sc Misha Verbitsky\\
{\sc Higher School of Economics\\
Faculty of Mathematics, 7 Vavilova Str. Moscow, Russia,}\\
\tt  verbit@mccme.ru}

}

\end{document}